\documentclass{sig-alternate}

\usepackage{multirow}

\usepackage[boxruled,boxed,lined]{algorithm2e}
\newcommand{\draft}{false}
\newcommand{\R}{\mathbb{R}}
\newcommand{\IR}{\mathbb{IR}}
\newcommand{\diff}{\backslash}

\newcommand{\ix}{{\bf x}}
\newcommand{\iy}{{\bf y}}
\newcommand{\imid}{\mathrm{mid}}
\newcommand{\iwid}{\mathrm{wid}}
\newcommand{\lb}[1]{\underline{#1}}
\newcommand{\ub}[1]{\overline{#1}}
\newcommand{\sol}{\mathrm{Sol}}
\newcommand{\cy}[1]{c^{#1}}
\newdef{example}{Example}
\newdef{remark}{Remark}
\newcommand{\bb}{{$2B$--consistency}}

\newcommand{\calx}[1]{\ensuremath{{\mathcal #1}}}

\begin{document}

\conferenceinfo{SAC'08}{March 16-20, 2008, Fortaleza, Cear\'{a}, Brazil}
\CopyrightYear{2008}
\crdata{978-1-59593-753-7/08/0003}

\title{An Efficient Algorithm for a Sharp Approximation of Universally Quantified Inequalities}

\numberofauthors{3}
\author{
\alignauthor A. Goldsztejn\\
       \affaddr{CNRS, University of Nantes}\\
       \affaddr{France}\\
       \email{alexandre@goldsztejn.com}
\alignauthor C. Michel\\
       \affaddr{University of Nice}\\
       \affaddr{France}\\
       \email{cpjm@essi.fr}
\alignauthor M. Rueher\\
       \affaddr{University of Nice}\\
       \affaddr{France}\\
       \email{rueher@essi.fr}
}

\date{14 September 2007}

\maketitle

\begin{abstract}
	This paper introduces a new algorithm for solving a sub-class of quantified constraint satisfaction problems (QCSP) where existential quantifiers precede universally quantified inequalities on continuous domains. This class of QCSPs has numerous applications in engineering and design. We propose here a new generic branch and prune algorithm for solving such continuous QCSPs. Standard pruning operators and solution identification operators are specialized for universally quantified inequalities. Special rules are also proposed for handling the parameters of the constraints. First experimentation show that our algorithm outperforms the state of the art methods.
\end{abstract}



\section{Introduction}

Many applications in engineering require verifying safety or performance conditions on a given system. For instance verifying for all electrical current in the interval $[i_\mathrm{min},i_\mathrm{max}]$ and all possible resistor values in the interval $r_0\pm5\%$ that the voltage across the resistor remains lower than a safety bound $u_\mathrm{max}$. In other words, we have to check the satisfiability of
\begin{equation}
	(\forall i\in[i_\mathrm{min},i_\mathrm{max}]) \ (\forall r\in[0.95\,r_0,1.05\,r_0]) \ (ri\leq u_\mathrm{max}).
\end{equation}	
Examples can be found in \cite{Benhamou2004} where the safety condition to be verified is a minimal mutual distance between some satellites during a complete revolution.

In design problems, requirements are a bit more general. Instead of verifying theses conditions, 
we have to determine values of design variables that satisfy safety and performance conditions 
for all values of uncertain physical data. 
All these problems can be modeled in the framework of quantified constraint satisfaction problems (QCSPs). 
QCSPs arise naturally in numerous other applications 
(cf. \cite{Jaulin1996,Jirstrand1997,Benhamou2004,Goldsztejn-SAC2006,Goldsztejn-Jaulin-CP2006,Ratschan}).

In this paper, we are interested in a restricted form of QCSPs where existential quantifiers precede universal quantifiers, i.e.
\begin{equation}\label{eq:QCSP}
	\exists x\in \ix, \forall y\in \iy, c_1(x,y)\wedge\cdots\wedge c_p(x,y),
\end{equation}
where $x=(x_1,\ldots,x_n)$ and $y=(y_1,\ldots,y_m)$ are vectors of variables, 
$\ix=(\ix_1,\ldots,\ix_n)$ and $\iy=(\iy_1,\ldots,\iy_m)$ are vectors of intervals over continuous domains, 
and constraints $c_i$  are inequalities of the form $f_i(x,y)\leq0$. 
In addition to the problems introduced in \cite{Benhamou2004}, this class of QCSPs can tackle 
the design of robust controllers, which has been addressed in numerous works 
\cite{Fiorio1993-CDC,Jaulin1996,Dorato1997,Malan1997,Garloff1998,Dorato2000}.

In general, we are not only interested in the decision problem (\ref{eq:QCSP}) 
but also in finding assignments of the existential variables which satisfy the constraints.

In other words, we consider $\ix$ as a search space where we try to find values $x\in\ix$ such that the relation $\forall y\in \iy, c_1(x,y)\wedge\cdots\wedge c_p(x,y)$ holds. Such a value $x$ is called a solution of (\ref{eq:QCSP}) and we define the solution set of (\ref{eq:QCSP}) in the following way:
\begin{equation}\label{eq:solution-set}
	\sol:=\{x\in \ix \ : \ \forall y\in \iy \ , \bigwedge_{i\in\{1,\ldots,p\}} c_i(x,y)\}.
\end{equation}

\begin{figure}
	\centering
	\epsfig{file=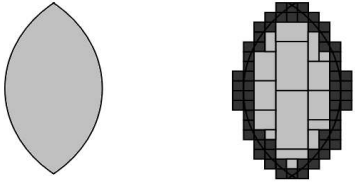, height=1in, width=2in, draft=\draft}
	\caption{\label{fig:paving}On the left hand side, the solution set $\sol$; on the right hand side, the inner approximation $\mathcal I$ in light gray and the boundary approximation $\mathcal B$ in dark gray.}
\end{figure}

An essential observation is that in all mentioned applications, solution sets are continuums of solutions
and that users are not really interested by isolated solutions but by continuous subsets of $\ix$ where all points are solutions. 
Indeed, in design problems determining the real value of a physical value without any tolerance information does not really make sense. 
That's why we compute two sets $\mathcal I$ and $\mathcal B$ (called respectively the \emph{inner approximation} and 
the \emph{boundary approximation}) which satisfy the following relation:
\begin{equation}
	\mathcal I \ \subseteq \ \sol \ \subseteq \ \mathcal I\cup\mathcal B.
\end{equation}
Set $\mathcal I$ contains only solutions whereas set $\mathcal B$ contains the boundary of the solution set. Thus $\mathcal B$ contains indifferently solutions and non-solutions. $\mathcal I\cup\mathcal B$ contains the whole solution set and is called the \emph{outer approximation}. Actually, due to computational limitations, $\mathcal I$ and $\mathcal B$ are unions of boxes (cf. Figure \ref{fig:paving}).

In this paper, we propose a new algorithm for computing the sets of boxes $\mathcal I$ and $\mathcal B$. In contrast with other methods dedicated to the resolution of QCSPs, we handle the QCSP (\ref{eq:QCSP}) as a non quantified CSP with quantified constraints: we introduce a CSP equivalent to (\ref{eq:QCSP}) formed of constraints $\cy{\iy}(x)$ on the variables $x$ (the domain parameter $\iy$ is explicitly given in exponent for clarity). These constraints are quantified inequalities: $\cy{\iy}(x)\equiv\forall y\in\iy,f(x,y)\leq0$; $y$ are called the parameters and $\iy$ their domains. We propose here a generic branch and prune algorithm dedicated to continuous CSP with parametric constraints. This algorithm is based on a specific implementation for universally quantified inequalities of the following techniques:

\indent\indent -- Pruning the non solution set;

\indent\indent -- Identifying the solution set;

\indent\indent -- Handling of parameter domains.

Other methods \cite{Benhamou2000,Ratschan:01a,Benhamou2004,Ratschan2006} have been proposed to compute approximations of the solution set of (\ref{eq:QCSP}). The new algorithm we propose here is both much simpler than the previously proposed algorithms, and also much more efficient. First experimentation show that our algorithm outperforms the methods introduced in \cite{Benhamou2000,Ratschan:01a,Benhamou2004,Ratschan2006}.

\paragraph*{Notations}

Boldface symbols denotes intervals, e.g. $\ix=[\lb{x},\ub{x}]=\{x\in\R:\lb{x}\leq x\leq\ub{x}\}$ 
where $\lb{x}$ and $\ub{x}$ belong to a finite subset of $\R$ (usually the floating point numbers, cf. \cite{Goldberg91}). 
The set of these intervals is denoted by $\IR$ and the set of $n$ dimensional interval vectors (also called boxes) by $\IR^n$. 
The width of an interval vector $\ix=(\ix_1,\ldots,\ix_n)$ is $\iwid(\ix)=\max_i|\lb{x}_i-\ub{x}_i|$ and 
its midpoint is $\imid(\ix)$. 
$\tilde x$ denotes a value of $\ix$.
For two boxes $\ix$ and $\iy$, the interval hull is denoted by $\ix\vee\iy$ 
and is the smallest box which contains both $\ix$ and $\iy$ 
(note the difference between the union $[-1,1]\cup[2,3]=\{x\in\R:-1\leq x\leq 1 \vee 2\leq x\leq 3\}$ which is disconnected, and the interval hull $[-1,1]\cup[2,3]=[-1,3]$). 
Also, their set difference is denoted by $\ix\diff\iy=\{x\in\ix:x\notin \iy\}$.

\paragraph*{Outline of the paper}

Section \ref{s:interval-analysis} recalls some basics on CSP with continuous domains. Section \ref{s:algorithm} describes the key features of the algorithm we propose. First experimental results are given in Section \ref{s:experimentations}.


\section{Basics on CSP with continuous domains}\label{s:interval-analysis}

To tackle CSP with continuous domains, a key issue is to be able to prove properties on continuums of real numbers. Interval analysis handles this problem in an efficient way using only computations on floating point numbers.

We recall here some basics which are required to understand the paper. More details can be found in \cite{Neum90,Kearfott1996-2}. 

An interval contractor for a constraint $c$ with $n$ variables is a function
\begin{equation}
	\mathrm{Contract}_c:\IR^n\longrightarrow \IR^n,
\end{equation}
that satisfies the following properties: (1) $\mathrm{Contract}_c(\ix )\subseteq\ix $; (2) $\forall x \ \bigl(x\in\ix  \wedge c(x)\bigr)\Longrightarrow x\in\mathrm{Contract}_c(\ix)$. Such contractors can be implemented using various techniques \cite{Neum90,Lhomme1993,Benhamou1994,Collavizza1999} for standard inequality constraints. Examples and experimentation presented in this paper use contractors based on the {\bb}. 

{\bb}~\cite{Cle87,
Lhomme1993} (also known as hull consistency)
states a local property on the
bounds of the domain of a variable at a single constraint level.
A constraint $c$ is $2B$-consistent if, for any variable $x_i$,
there exist values in the
domains of all other variables that satisfy $c$ when $x_i$ is fixed
to $\underline{x}_i$ and $\overline{x}_i$.

The filtering by {\bb} of
$P= ({\bf x},\calx{C})$ is the CSP $P'= ({\bf x'},\calx{C})$ such that

\begin{itemize}
\item $P$ and $P'$ have the same solutions;

\item $P'$ is $2B$-consistent;

\item ${\bf x'} \subseteq {\bf x}$ and the domains in ${\bf x'}$ are
the largest ones for which $P'$ is $2B$-consistent.
\end{itemize}

\noindent
Filtering by {\bb} of $P$ always exists and is unique~\cite{Lhomme1993},
that is to say, it is a closure.\\


\section{Description of the Algorithm}\label{s:algorithm}

The key idea of our algorithm is to reformulate a QCSP as a CSP with parametric constraints.
More precisely, we reformulate a QCSP with constraints of type (\ref{eq:QCSP}) 
as a CSP $\langle \mathcal V,\mathcal C,\mathcal D\rangle$ where:

\indent\indent -- $\mathcal V=(x_1,\ldots,x_n)$;

\indent\indent -- $\mathcal C=\{\cy{\iy}_{1},\ldots,\cy{\iy}_p\}$ with $\cy{\iy}_{i}(x)\equiv\forall y\in\iy,f_i(x,y)\leq0$;

\indent\indent -- $\mathcal D=(\ix_1,\ldots,\ix_n)$.

The parameter domains are considered at the level of each constraint, 
so we can process them differently for each constraint; they are initialized with $\iy$. 
Standard techniques have to be adapted for this class of constraints. 
For example, we show in Subsection \ref{ss:pruning} that
$\mathrm{Contract}_{\forall y\in\iy,f(x,y)\leq0}(\ix)$ can be achieved with $\mathrm{Contract}_{f(x,\tilde y)\leq0}(\ix)$
for any $\tilde y\in\iy$, the latter being implemented using standard interval contractors.

\subsection{Outline of the Algorithm}

The branch and prune algorithm we propose alternates pruning and branching steps in a standard way to reject parts of the search space that do not contain any solution. This process is interleaved with the identification of inner boxes that contain only solutions. 

We also introduce new parameter instantiations and parameter reduction techniques
which play a key role during the pruning of the search space and
the identification of the solution sets.
These key points are 
detailed in subsections \ref{ss:pruning},  \ref{ss:solution-identification} and \ref{ss:parameter-reduction}.
They are implemented within a branching process (see Algorithm \ref{alg:bb}) using the following functions:

\incmargin{1em}
\linesnumbered
\begin{algorithm}
	\KwIn{$\mathcal C$, $\ix $, $\epsilon$}
	\KwOut{$(\mathcal I,\mathcal B)$}
	$\mathcal{U}\;\leftarrow \{(\mathcal C,\ix )\}$;
	$\mathcal{B}\;\leftarrow \{\}$;
	$\mathcal{I}\;\leftarrow \{\}$\;
	\While{$\neg\mathrm{empty}(\mathcal{U})$}
	{
		$(\mathcal C,\ix)\leftarrow\mathrm{extract}(\mathcal{U})$\;
		\uIf{$\iwid(\ix)>\epsilon$}
		{
			$\mathcal C'\leftarrow\mathrm{ParameterInstantiation}(\mathcal C,\ix)$\;
			$\ix'\leftarrow\mathrm{Pruning}(\mathcal C',\ix)$\;
			$(\mathcal C'',\ix'',\mathcal I')\leftarrow\mathrm{SolutionIndentification}(\mathcal C',\ix')$\;
			$\mathcal I\leftarrow\mathcal I\cup\mathcal I'$\;
			$\mathcal C'''\leftarrow\mathrm{ParameterDomainBisection}(\mathcal C'',\ix'')$\;
			$\mathcal U\leftarrow\mathcal U\cup\mathrm{Branching}(\mathcal C''',\ix'')$\;
		}
		\Else
		{
			$\mathcal{B}\leftarrow \mathcal{B}\cup\{\ix \}$\;
		}
	}
	\Return{$(\mathcal{I},\mathcal{B})$}\;
	\caption{\label{alg:bb}Generic Branch and Prune Algorithm}
\end{algorithm}

Function \emph{ParameterInstantiation()} returns a set of constraints $\mathcal C'$ which is equivalent to $\mathcal C$ on the domain $\ix$, but where some parameters have been instantiated (see Subsubsection \ref{sss:pdi}).

Function \emph{Pruning()} contracts the domain $\ix$ without loosing any solution of $\mathcal C'$ (see Subsection \ref{ss:pruning}).

Function \emph{SolutionIndentification()} contracts $\ix'$ to $\ix''$ and returns a list of inner boxes $\mathcal I'$ (see Subsection \ref{ss:solution-identification}). The domains of the parameters are also updated in $\mathcal C''$ (see Subsubsection \ref{sss:pdp}).

Function \emph{ParameterDomainBisection()} bisects the parameter domains of the constraints and thus increases the number of constraints in $\mathcal C''$ (see Subsubsection \ref{sss:pdb}).

Function \emph{Branching()} achieves a standard bisection of the domains $\ix''$.

\subsection{Pruning}\label{ss:pruning}

\paragraph*{Local Pruning}

Given a constraint $\forall y\in\iy^i,c_i(x,y)$ and a box $\ix$, we want to contract $\ix$ rejecting only parts which do not contain any solution of this constraint. To achieve this task, Benhamou et al. and Ratschan \cite{Benhamou2004,Ratschan2006} apply $\mathrm{Contract}_{c_i(x,\iy^i)}(\ix)$, where $y$ is handled as an existentially quantified variable in the domain $\iy^i$. However, this strategy lacks efficiency. That's why we propose a better contractor by instantiating the parameter to an arbitrary value $\tilde y\in\iy^i$ (see Example \ref{ex:outer}). More formally, the box $\ix$ is contracted using $\mathrm{Contract}_{c_i(x,\tilde y)}(\ix)$ which rejects only parts of $\ix$ that do not satisfy $c_i(x,\tilde y)$, and thus, do not satisfy $\forall y\in\iy^i,c_i(x,y)$.

Many strategies can be used to choose $\tilde y$. We chose $\tilde y=\imid(\iy^i)$. Experimentation have shown that searching for a better value of $\tilde y$ is not worthwhile. This is due to the fact that the parameter handling techniques reduce the impact of this choice (see examples \ref{ex:pdb-outer-contraction} and \ref{ex:dbpdr-outer-contraction} in Subsection \ref{ss:parameter-reduction}). 

\begin{example}\label{ex:outer}
	Let us consider the constraint $c(x)$ defined by $\forall y\in\iy,f(x,y)\leq0$ with $f(x,y)=10y-x-y^2$, $\iy=[0,1]$, and $\ix=[0,15]$. The solution set of this simple CSP is the interval $[9,15]$. To reject values of $\ix$ that do not satisfy $c(x)$, we apply $\mathrm{Contract}_{f(x,0.5)\leq0}(\ix)$, which reduces $\ix$ to $\ix'=[4.75,15]$. The method proposed in \cite{Benhamou2004,Ratschan2006} computes $\mathrm{Contract}_{f(x,\iy)\leq0}(\ix)$ but cannot achieve any contraction.
\end{example}

\paragraph*{Global Pruning}

\begin{figure*}
	\centering
	\epsfig{file=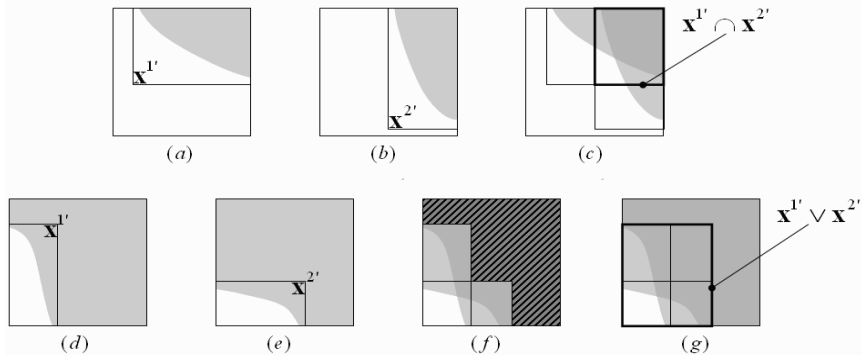, height=1.9in, width=4.5in, draft=\draft}
	\caption{\label{fig:global-prunings}Upper row: pruning computed from local contractions. Lower row: identification of solutions using local contractions of the negation of the constraints.}
\end{figure*}

The contraction of $\ix$ using one constraint removes solutions of this constraint, and therefore solutions of the global CSP. This is illustrated in the first row of Figure \ref{fig:global-prunings} where diagrams (a) and (b) show two local contractions leading to $\ix^{1\prime}$ and $\ix^{2\prime}$. The global contraction depicted in Diagram (c) is obtained by computing $\ix^{1\prime}\cap\ix^{2\prime}$. Of course, in practice we compute this intersection in an incremental way.

\subsection{Identification of Sets of Solutions}\label{ss:solution-identification}

As in \cite{Collavizza1999-b,Ratschan:01a,Vu2002,Benhamou2004,Ratschan2006,Vu2006}, the identification of sets of solutions is implemented by applying interval contractors to the negation of the constraints. Again, this process is achieved for each constraint separately. However, now only areas that are proved to be local solutions for all constraints are solutions of the whole CSP.

\paragraph*{Identification of Solutions for a Single Constraint}

In order to identify parts of a box $\ix$ which contains only solutions of the constraint $\forall y\in\iy^i,f_i(x,y)\leq 0$, we compute
\begin{equation}\label{eq:inner-pruning}
	\langle\ix^{i\prime},\iy^{i\prime}\rangle=\mathrm{Contract}_{f_i(x,y)\geq0}(\langle\ix,\iy^i\rangle),
\end{equation}
where $\langle\ix,\iy\rangle$ stands for the vector $(\ix_1,\ldots,\ix_n,\iy_1,\ldots,\iy_m)$. Thus, every value of $\ix\diff\ix^{i\prime}$ satisfies $f_i(x,y)<0$ for all values of $y$ in $\iy^i$. Indeed, the part of $\langle \ix,\iy^i \rangle$ which has been pruned does not contain any solution:
\begin{equation}\label{eq:inner-pruning-2}
	\forall x\in\ix \ \forall y\in\iy^i \ , \ (x\notin\ix^{i\prime}\vee y\notin\iy^{i\prime})\Longrightarrow \neg (f_i(x,y)\geq0).
\end{equation}
From (\ref{eq:inner-pruning-2}) we have trivially
\begin{equation}
	\forall x\in\ix \ \forall y\in\iy^i \ , \ x\notin\ix^{i\prime} \ \Longrightarrow \ f_i(x,y)<0.
\end{equation}
This process is illustrated on the following example.
\begin{example}\label{ex:inner}
	Continuing Example \ref{ex:outer}, let us consider $c(x)$ and $\ix=[4.75,15]$. In order to identify values of $\ix$ that satisfy $c(x)$, we apply $\mathrm{Contract}_{f(x,y)\geq0}(\ix,\iy)$ and obtain $\ix'=[4.75,10]$ and $\iy'=[0.475,1]$. Therefore, $\ix\diff\ix'=]10,15]$ contains only solutions of the constraint.
\end{example}

\paragraph*{Identification of Solutions of the whole CSP}

The $\ix^{i\prime}$ for $i\in\{1,\ldots,p\}$ have been computed in such a way that all values of $\ix\diff\ix^{i\prime}$ 
verify $c_i(x)$ (cf. the diagrams (d) and (e) of the second row of Figure \ref{fig:global-prunings}). 
Thus, the values of $\ix$ that are outside all $\ix^{i\prime}$ satisfy all the constraints. Formally,
\begin{equation}\label{eq:inner-union}
	\ix\diff(\ix^{1\prime}\cup\cdots\cup\ix^{p\prime})
\end{equation}
contains only solutions of the CSP. This inner approximation (\ref{eq:inner-union}) is the dashed area of Diagram (f) of the second row of Figure \ref{fig:global-prunings}. Its description becomes very complicated for higher dimensions and when numerous constraints are involved. As in \cite{Ratschan2006}, we use instead the weaker inner approximation
\begin{equation}\label{eq:inner-join}
	\ix\diff(\ix^{1\prime}\vee\cdots\vee\ix^{p\prime}).
\end{equation}
This inner approximation 
is represented on Diagram (g) of the second row of Figure \ref{fig:global-prunings}. Let us define $\ix'=\ix^{1\prime}\vee\cdots\vee\ix^{p\prime}$. The closure\footnote{The set difference $\ix\diff\ix'$ is not closed in general, e.g. $[-10,10]\diff[-1,1]=[-10,-1[\cup]1,10]$. Nevertheless, we can observe that if $f$ is continuous and non positive in $\ix\diff\ix'$, then it is also non positive in its closure $\mathrm{cl}(\ix\diff\ix')$, e.g. if $f(x)\leq0$ in $[-10,-1[\cup]1,10]$ then $f(x)\leq0$ also holds in $[-10,-1]\cup[1,10]$. Thus, $\mathrm{cl}(\ix\diff\ix')$, which can be described by a simple union of boxes, is recorded instead of $\ix\diff\ix'$.} of $\ix\diff\ix'$ is added to the set of inner boxes while $\ix'$ has to be further explored.

\subsection{Handling Parameter Domains}\label{ss:parameter-reduction}

As said before, the size of the domains of the parameters is a critical issue for an efficient application of interval contractors. This section details the three methods implemented in our algorithm to overcome this problem. The goal of these methods is to only apply interval contractors to parts of the initial parameter domains while keeping the CSP solution set unchanged.

\subsubsection{Parameter Domain Pruning}\label{sss:pdp}

During the identification of solutions (\ref{eq:inner-pruning}) the domains of the parameters $\iy^i$ are reduced to $\iy^{i\prime}$. 
In \cite{Benhamou2004,Ratschan2006}, the initial parameter domain is restored, so this contraction of the parameter
domains is not propagated. However, this reduced domain can be used while keeping the solution set of the constraint unchanged. 
Indeed, (\ref{eq:inner-pruning-2}) trivially entails
\begin{equation}
	\forall x\in\ix \ \forall y\in(\iy^i\diff\iy^{i\prime}) \ f_i(x,y)<0,
\end{equation}
As a consequence, for every fixed $x$ in $\ix$, $\forall y\in\iy^{i\prime},c(x,y)$ implies $\forall y\in\iy^i,c(x,y)$. Finally, we can replace the latter by the former keeping the solution set unchanged during the solution identification step of the algorithm. This is a significant improvement that is implemented without any overhead since the parameter domains have to be contracted anyway during the solution identification step. 

\begin{example}
	In Example \ref{ex:inner}, the identification of solutions step has contracted the parameter domain from $\iy=[0,1]$ to $\iy'=[0.475,1]$. This latter domain can be used in future processing while keeping unchanged the CSP solution set.
\end{example}

\subsubsection{Parameter Domain Bisection}\label{sss:pdb}

Bisecting parameter domains is the most obvious way to reduce them. This is a critical issue to improve the convergence of the algorithm. Bisection is also used in \cite{Ratschan2006} but implemented in a different way. Here, we change the constraint $\forall y\in\iy,c(x,y)$ to the conjunction of the two constraints $\forall y\in\iy^1,c(x,y)$ and $\forall y\in\iy^2,c(x,y)$, where $\iy^1$ and $\iy^2$ are obtained by bisecting $\iy$. The following example illustrates how the parameter domain bisection is performed, and how it improves the pruning process as well as the identification of solutions.

\begin{example}\label{ex:pdb-outer-contraction}
	Let us come back to Example \ref{ex:outer} where the constraint store contains only one constraint $\mathcal C=\{(\forall y\in\iy,f(x,y)\leq0)\}$. Bisecting parameter domains, we obtain the new constraint store $\mathcal C'=\{(\forall y\in\iy^1,f(x,y)\leq0), (\forall y\in\iy^2,f(x,y)\leq0)\}$ with $\iy^1=[0,0.5]$ and $\iy^2=[0.5,1]$. The latter constraint store is equivalent to the original one, but contains more constraints with smaller parameter domains.
	
	The pruning operator is now applied to each constraint of the store, thus the two contractions $\mathrm{Contract}_{f(x,0.25)\leq0}(\ix)$ and $\mathrm{Contract}_{f(x,0.75)\leq0}(\ix)$ are computed. This leads to $\ix'=[6.9375,15]$. The pruning achieved here is sharper than the one computed without bisecting the parameter domains (cf. Example \ref{ex:outer}).
	
	To identify solutions of this new CSP, we compute $(\ix^{i\prime},\iy^{i\prime})=\mathrm{Contract}_{f(x,y)\geq0}(\ix',\iy^i)$ for $i\in\{1,2\}$. We obtain $\ix^{1\prime}=\emptyset$, $\iy^{1\prime}=\emptyset$, $\ix^{2\prime}=[6.9375,9.75]$ and $\iy^{2\prime}=[0.71875,1]$. Finally, $\ix'$ is contracted to $\ix^{1\prime}\vee\ix^{2\prime}=[6.9375,9.75]$ while $\mathrm{cl}([6.9375,15]\diff[6.9375,9.75])=[9.75,15]$ is proved to be an inner box. The solution identification is also more efficient thanks to the bisection of the parameter domain (cf. Example \ref{ex:inner}).
\end{example}

\subsubsection{Parameter Instantiation}\label{sss:pdi}

The constraint $\forall y\in\iy,f(x,y)\leq0$, where $\iy=[\lb{y},\ub{y}]$, can be simplified if the function $f$ is proved to be monotonic w.r.t. the parameter. Indeed, if $f$ is increasing (resp. decreasing) w.r.t. $y$, then the quantified constraint is obviously equivalent to the non-quantified constraint $f(x,\ub{y})\leq0$ (resp. $f(x,\lb{y})\leq0$). The variation of the function can be checked by evaluating its derivative w.r.t. the parameter on the intervals $\ix$ and $\iy$. The same property holds if there are several parameters, using a derivative w.r.t. each parameter.

\begin{example}\label{ex:dbpdr-outer-contraction}
	Let us come back to Example \ref{ex:outer}. Evaluating
	\begin{equation}
		(\partial f/\partial y)(\ix,\iy)=10-2\iy=[8,10],
	\end{equation}
	we prove that $f$ is increasing w.r.t. $y$. Therefore, the constraint $\forall y\in\iy,f(x,y)\leq0$ is equivalent to $f(x,1)\leq0$. Then, the pruning contracts $\ix=[0,15]$ to $[9,15]$ and the solution identification proves that $[9,15]$ contains only solutions. On this simple example, we obtain an exact description of the solution set. This set is much sharper than the contractions obtained without this instantiation of the parameter (cf. examples \ref{ex:outer} and \ref{ex:inner}).
\end{example}

The parameter instantiation can drastically improve the efficiency of the pruning and the solution identification steps. However, evaluating derivatives can be expensive, in particular when they involve trigonometric function. Experimentation presented in the next section illustrate well this trade-off.


\section{Experimentation}\label{s:experimentations}

\begin{center}
\begin{table}[t]
\begin{center}
\small{\begin{tabular}{|r|r||r||r|r|}
\hline
 & & & \multicolumn{2}{c|}{Qine} \\ \hline
problem & ratio & Rsolver & 2B & 2B+\\ \hline\hline

\multirow{3}{*}{Circle} & 0.98 & $55.67$ & $0.90$ & $1.03$ \\
  & 0.99 & $169.51$ & $2.32$ & $2.22$ \\
  & 0.999 & $-$ & $69.18$ & $31.55$ \\ \hline

\multirow{3}{*}{PathPoint} & 0.6 & $122.72$ & $0.01$ & $0.01$ \\
  & 0.65 & $334.34$ & $0.01$ & $0.02$ \\
  & 0.7 & $-$ & $0.02$ & $0.03$ \\ \hline

\multirow{3}{*}{Parabola} & 0.96 & $48.33$ & $1.47$ & $1.36$ \\ 
  & 0.97 & $130.30$ & $2.74$ & $1.92$ \\
  & 0.98 & $-$ & $8.77$ & $6.59$ \\ \hline

\multirow{3}{*}{Robot} & 0.98 & $10.34$ & $0.06$ & $0.08$ \\
  & 0.99 & $26.74$ & $0.14$ & $0.18$  \\
  & 0.999 & $-$ & $5.37$ & $4.08$  \\ \hline

\multirow{3}{*}{Satellite} & 0.5 & $303.30$ & $71.36$ & $114.33$ \\
  & 0.55 & $-$ & $168.32$ & $268.36$ \\
  & 0.6 & $-$  & $227.09$ & $368.90$ \\ \hline

\multirow{3}{*}{Robust1} & 0.999 & $1.10$ & $0.01$ & $0.00$ \\
  & 0.9999 & $7.16$ & $0.04$ & $0.00$  \\
  & 0.99999 & $76.25$ & $0.10$ & $0.00$ \\ \hline

\multirow{3}{*}{Robust4} & 0.5 & $-$ & $-^*$ & $0.00$ \\
  & 0.55 & $-$ & $-$ & $0.01$ \\
  & 0.6 & $-$ & $-$  & $0.01$ \\ \hline

\multirow{3}{*}{Robust5} & 0.7 & $75.35$ & $0.00$ & $0.00$ \\
  & 0.75 & $80.62$ & $0.00$ & $0.00$ \\
  & 0.8 & $-$ & $0.01$ & $0.00$ \\ \hline

\multirow{3}{*}{Robust6} & 0.7 & $59.98$ & $0.02$ & $0.00$ \\
  & 0.75 & $224.33$ & $0.05$ & $0.00$  \\
  & 0.8 & $-$ & $0.09$ & $0.00$ \\ \hline
\end{tabular}}
\end{center}
\caption{Timing (in seconds) for Rsolve vs Qine}\label{benches}
\end{table}
\end{center}

This section compares the results of the IPA system described in \cite{Benhamou2004},  
Rsolver \cite{Ratschan2006} from S. Ratschan with our system, called ``Qine'',
on a set of $9$ benches.

The $9$ benches come from the literature.
The Circle, PathPoint, Parabola, Robot and Satellite QCSPs are taken from \cite{Benhamou2004}.
A description of Robust1 (respectively, Robust4, Robust5 and Robust6) can
be found in \cite{Dorato2000} (respectively, \cite{Malan1997}, \cite{Jaulin1996} and \cite{Jaulin2002-CDC}).

The Rsolver and Qine benches have been run on an Intel Core Duo 2 at 2.4Ghz with a time out of 600s.
To limit the effect of the high memory consumption of these algorithms,
the available memory has been restricted to 1Gb.
Thus, a bench could either succeed to run within these two limits,
end with a time out (``$-$''), or reach the memory limit (``$-^*$'').

Table \ref{benches} reports the results obtained with Rsolver and Qine.
It gives Rsolver timing, as well as the time required to solve the benches
for the different Qine running options :
\begin{itemize}
\item ``2B'' uses the contractor based on $2b$-consistency techniques.
More precisely, its implementation relies on a forward-backward evaluation
of the direct acyclic graph which represents the constraint.
\item ``2B+'' combines the previous contractor with 
the derivative based parameter handling strategy introduced in subsection \ref{sss:pdi}.
\end{itemize}
The results have been computed according to a ratio (column 2) where
\begin{equation}
ratio = \frac{V_{inner} + V_{outer}}{V_{initial}}
\end{equation}
where $V_{inner}$ is the total volume of the inner boxes, 
 $V_{outer}$ is the total volume of the outer boxes and
 $V_{initial}$ is the volume given by the initial domains of the variables.
Though all the benches have been run for ratios going from $0.5$ to $0.99999$,
the table gives only the most significant results for the sake of space.

As shown in table \ref{benches}, in average, Qine outperforms Rsolver by one order of magnitude.
For instance, Qine handles the Robust benches immediately while Rsolver needs much for
time to do so. 
A comparison of the different available combinations shows
that the ``2B+'' combination has a more robust behavior.
It takes advantage of all the available information and offers a good
trade-off between the computation time and the domain reductions.
However, on the Satellite bench, 2B performs better than 2B+.
This illustrates the trade-off between the cost of derivative computation
and the benefit of parameter instantiations.

\begin{center}
\begin{table}[t]
\begin{center}
\small{\begin{tabular}{|l|r||r|r||r|r|}
\hline
 & & \multicolumn{2}{c||}{IPA} & \multicolumn{2}{c|}{QINE}\\ \hline
problem & ratio & $t$ & $t*corr$ & 2B & 2B+\\ \hline\hline
\multirow{3}{*}{Circle} 
  & 0.8928 & 2.4 & 1.402 & 0.14 & 0.20\\ 
  & 0.9326 & 9.328 & 5.450 & 0.20 & 0.30\\ 
  & 0.9535 & 64.876 & 37.907 & 0.30 & 0.42\\ 
\hline
PathPoint 
  & 0.8172 & 148.66 & 86.86 & 0.08 & 0.12\\ 
\hline
\multirow{3}{*}{Parabola} 
  & 0.8716 & 0.340 & 0.1986 & 0.06 & 0.06\\ 
  & 0.9340 & 2.824 & 1.65 & 0.34 & 0.33\\ 
  & 0.9650 & 75.936 & 44.3697 & 1.80 & 1.65\\ 
\hline
\multirow{3}{*}{Robot} 
  & 0.9924 & 1.112 & 0.6497 & 0.26 & 0.34\\ 
  & 0.9973 & 5.908 & 3.452 & 1.07 & 1.10\\ 
  & 0.9980 & 16.933 & 9.894 & 1.89 & 1.81\\ 
\hline
\multirow{3}{*}{Satellite} 
  & 0.2813 & 5.660 & 3.3071 & 18.41 & 29.61\\ 
  & 0.4844 & 27.633 & 16.146 & 52.79 & 84.90\\ 
\hline
\end{tabular}}
\end{center}
\caption{Timing (in seconds) for IPA vs Qine}\label{benches2}
\end{table}
\end{center}

Tables \ref{benches2} compares IPA with Qine. 
IPA system has been run on a Pentium M at 2Ghz running Linux.
To allow a fair comparison, we have computed a timing correction :
the same system, Rsolver, has been run on both systems
in order to determine this correction.
Therefore, the initial timing obtained for IPA (column 3) has been
multiplied by 0.58 (column 4) to allow a fair comparison
between the two systems.

Here again, Qine outperforms IPA. Indeed, IPA was not
able to solve the PathPoint bench within the timeout for only one of the tested ratios.
However, IPA is faster than Qine on the Satellite bench for the
lower ratio values. This behavior is
probably due to the limit of the 2B based contractor whose domain
reduction capabilities decrease when a variable
has multiple occurrences within one constraint.
IPA is based on a Box contractor which does not
suffer from the same behavior (see \cite{Collavizza1999} for
a detailed comparison of 2B and Box).
However, when the ratio value increases, 
Qine becomes faster than IPA. For instance,
IPA needs more than 876s (corrected time) to solve
the Satellite bench for a ratio of 0.6276
whereas Qine achieves this task with 2B within 401.44s.

Note that the class of QCSPs handled by IPA is limited to one parameter.
Thus, IPA is not able to solve the Robust benches.


\section{Conclusion}

In this paper, we have introduced a new, simple and efficient algorithm
to handle a significant class of QCSPs. 
Examples coming from the literature reveal that this
class covers most of the practical applications.

Our algorithm is based on new techniques for handling parameters.
It also takes advantage of information provided by the derivatives
to improve the contraction of the domains involved in constraints with parameters,
a key issue in the efficient solving of QCSPs.

Experimentation underline the efficiency of our algorithm which outperforms
two of the available state of the art implementations able to handle such QCSPs.
In average, our implementation is by one order of magnitude quicker than the
two other systems.

Further work concerns the improvement of the implementation of the contractor
by using the best filtering techniques for each type of constraints
and the generalization of
the use of available information to still enhance the speed of the
solving process.

\section*{Acknowledgments} 
We are grateful to Marc Christie for his valuable help in the experimentation.


\bibliographystyle{abbrv}

\end{document}